\documentclass[a4paper,10pt]{article}
\usepackage[cp1251]{inputenc}
\usepackage[english]{babel}

\usepackage{amssymb}\usepackage{amsmath}
\usepackage{epsfig}
\def\bl{\rule[-1mm]{2.4mm}{2.4mm}}

\def\be{\begin{equation}}
\def\ee{\end{equation}}
\newtheorem{theorem}{\bf Theorem}
\newtheorem{lemma}{\bf Lemma}

\begin{document}

\title {Rational functions admitting double decomposition} 
\author{\copyright 2010 ~~~~A.B.Bogatyrev
\thanks{Partially supported by RFBR grants 10-01-00407 and RAS Program 
"Modern problems of theoretical mathematics"}} 
\date{} 
\maketitle

J.Ritt \cite{R1} has investigated the structure of complex polynomials with respect 
to superposition. The polynomial $P(x)$ is said to be indecomposable iff the representation
$P=P_1\circ P_2$ means that either $P_1$ or $P_2$ is a linear function. The decomposition 
$P=P_1\circ P_2\circ \dots\circ P_r$ is called maximal if all factors $P_j$ are indecomposable polynomials
and are not linear.  Ritt proves that any two maximal decompositions have the same length $r$, the same 
(unordered) set $\{\deg (P_j)\}$ of factor's degrees and may be connected by a finite chain of transformations, 
each step consists in replacing the left side of the following double decomposition
\be 
R_1\circ R_2=R_3\circ R_4
\label{DD1}
\ee 
by its right side. The solutions of the latter functional equation are indecomposable polynomials of degrees greater 
than one and all of them were explicitly listed by Ritt. 

The analogues of Ritt theory for rational functions were constructed just for several particular 
classes of the said functions, say for Laurent polynomials \cite {FP1}. In this note we describe a certain class of 
double decompositions (\ref{DD1}) with rational functions $R_j(x)$ of degree greater than one.
Essentially, described below rational functions were discovered by E.I.Zolotarev in 1877 
as a solution of certain optimization problem \cite{Z1, A1}. However, the double decomposition property for them  
was hidden until recently because of somewhat awkward representation. Below we give a (possibly new) symmetric representation of Zolotarev fractions resembling the parametric representation for Chebyshev polynomials, which are a special limit case of Zolotarev fraction.

\section{Zolotarev fractions and their nesting property}
Let $\tau\in i\mathbb{R}_+$ and $\Pi(\tau)$ be a rectangle of size $2\times|\tau|$:
$$
\Pi(\tau):=\{u\in\mathbb{C}:\quad |Re~u|\le1, 0\le Im~u\le|\tau| \}.
$$
The conformal mapping $x_\tau(u)$ of this rectangle to the upper half plane fixing three points 
$u=\pm 1, 0,$ has a very simple appearance 
$$
x_\tau(u)=sn(K(\tau)u|\tau)
$$
in terms of elliptic sine $sn$ and complete elliptic integral $K$. From the reflection principle
for conformal mappings it may be easily derived that the parametric representation:
$$
R(u):=x_\tau(u);\quad x(u):=x_{n\tau}(u), \qquad u\in\mathbb{C}, \qquad n\in\mathbb{N},
$$
gives a degree $n$ rational function $R$ of argument $x$:
$$
Z_n(x|\tau):=R(u(x))=x_\tau\circ x_{n\tau}^{-1}.
$$
This rational function is known as  Zolotarev fraction. 
Directly from the definition it follows that Zolotarev fractions obey the nesting property:
\be
Z_{mn}(x|\tau)=Z_m(Z_n(x|m\tau)|\tau), \qquad m,n\in \mathbb{N}. 
\label{NPr}
\ee
When parameter $\tau$ tends to zero (suitably renormalized) Zolotarev fraction becomes classical Chebyshev polynomial 
and the well known nesting property of Chebyshev polynomials becomes just the consequence of the above formula.
Interchanging $n$ and $m$ in formula (\ref{NPr}) we observe that Zolotarev fractions of composite degrees possess double decompositions of the kind (\ref{DD1}). We generalize the construction of Zolotarev fraction in the next section.

\section{Construction}
Let $ L$ be a rank two lattice in the complex plane of variable $u$. The group of translations of the plane by the elements of the lattice we designate by the same letter $ L$. Let $ L^+$ be the group $ L$ extended by degree two transformation $u\to -u$. The extended group acts discontinuously in the complex plane, so the orbit space is well defined and carries natural complex structure
$$
{\mathbb C}/ L^+={\mathbb C}P^1.
$$
We can introduce a global coordinate on this Riemann sphere, say 
$$
x(u)=\wp(u| L):=u^{-2}+\sum_{0\neq v\in L}((u-v)^{-2}-v^{-2}).
\label{coord}
$$
Some basis in the lattice $ L$ is traditionally used as the second argument of the Weierstrass function, 
however it depends on the lattice as a whole. 

Once we have a full rank sublattice $ L_\bullet$ of $ L$, the group $ L_\bullet^+$ is a subgroup of $ L^+$
and any orbit of $ L_\bullet^+$ is contained in the orbit of $ L^+$. Therefore we have a holomorphic mapping
of one sphere to the other:
\be
{\mathbb C}/ L_\bullet^+\to {\mathbb C}/ L^+,
\ee 
which becomes a rational function once we fix complex coordinate on each sphere. Thus we obtain a 
degree $| L: L_\bullet|$ rational function $R_{ L: L_\bullet}(x)$:
\be
R_{ L: L_\bullet}(x_\bullet(u)):= x(u),   
\qquad    x_\bullet(u):= \wp(u| L_\bullet),
\label{R}
\ee 
which is a general form of $g=1$ rational functions in the terminology of \cite{B1}.
To get modulus  $\tau\in i\mathbb{R}_+$ Zolotarev fraction  we just need to take $ L=Span_\mathbb{Z}\{4, 2\tau\}$ and
$ L_\bullet=Span_\mathbb{Z}\{4, 2n\tau\}$,  then $R_{ L: L_\bullet}(x)$ coinsides with  
$Z_n(x|\tau)$ up to normalization (i.e. pre- and post- compositions with linear fractional functions).

Suppose we have two different sublattices $ L_\bullet$ and $ L_\circ$ of the same lattice $ L$.
Their intersection $ L_{\bullet\circ}:= L_\bullet\cap L_\circ$ is a full rank sublattice of both $ L_\bullet$ and $ L_\circ$. Indeed, $ L_{\bullet\circ}$ contains a full rank sublattice $| L: L_\bullet|| L: L_\circ|~ L$.
Obviously, we have a double decomposition:
\be
\label{DD}
R_{ L: L_{\bullet\circ}}=R_{ L: L_\bullet}\circ R_{ L_\bullet: L_{\bullet\circ}}=
R_{ L: L_\circ}\circ R_{ L_\circ: L_{\bullet\circ}}.
\ee

Not all of the relations (\ref{DD}) are independent.
Below we show that arbitrary double decomposition (\ref{DD}) is a consequence of the 
same relations for prime index sublattices 
$ L_\bullet$, $ L_\circ$ of $ L$. 

\section{Prime index sublattices}
Given a base in the lattice $ L$, a base in its sublattice $ L_\bullet$ is obtained via two by two matrix $Q$ with integer entries. Other choice of bases results in multiplication of $Q$ by invertible integer matrices (i.e. of determinant $\pm1$) on the left and on the right. The index of sublattice $ L_\bullet$ in $ L$ denoted by $| L: L_\bullet|$ equals to $|det ~Q|$ and  it is independent of the choice of bases in the lattice and its sublattice. Given a chain of lattices $ L\supset  L_\bullet\supset  L_{\bullet\bullet}$, the indecies obey the multiplication rule:  $| L: L_{\bullet\bullet}|=| L: L_\bullet|| L_\bullet: L_{\bullet\bullet}|$.

\begin{lemma}
Any prime index $p$ sublattice $ L_\bullet$ of $ L$ has the following representation
\be
\label{repL}
 L_\bullet=Span_{\mathbb Z}\{p L, e\}
\ee
where $e$ is any element of $ L_\bullet\setminus{p L}$. Conversely, the right hand side of (\ref{repL})
is an index $p$ sublattice of $ L$ provided $e\not\in p L$.
\end{lemma}
{\it Proof.}
Let the matrix $Q\in GL_2(\mathbb{Z})$ maps the base of $ L$ to the base of $ L_\bullet$. The matrix $pQ^{-1}$ is integer and therefore $ L_\bullet$ contains sublattice $p L$ of the same index $p$. We get the following chain of sublattices
$$
p L\subset Span_{\mathbb Z}\{p L, e\} \subset  L_\bullet
$$
Prime index $p=| L_\bullet:p L|$ is the product of indecies $| L_\bullet:Span\{\dots\}|$ and $|Span\{\dots\}:p L|$, therefore one of them should be unity. In other words, the middle lattice in the chain is equal either to the left or to the right lattice in the chain. The choice of the element $e$ says that the middle lattice in the chain is strictly larger than $p L$. ~~\bl

{\bf Corollary 1.} {\it Let $L_\bullet\neq L_\circ$ be two sublattices of $ L$ of the same prime index $p$. Then $L_\bullet\cap L_\circ=$ $p L$.}

{\it Proof}. Each index $p$ sublattice of $ L$ contains $p L$. If there is at least one more element $e$ in the intersection 
$L_\bullet\cap L_\circ$ then each of two sublattices may be reconstructed by formula (\ref{repL}) and therefore they coinside. ~~\bl

{\bf Corollary 2.} {\it Let $L_\bullet$ and $L_\circ$ be two sublattices of $ L$ of different prime indecies $p_\bullet$ and $p_\circ$ respectively. Then their intersection  has the representation:

\be
\label{repL12}
L_\bullet\cap L_\circ= Span_{\mathbb Z}\{p_\bullet p_\circ L, ~p_\bullet e_\circ, ~p_\circ e_\bullet\}
\ee
where $e_*$ is any element of $L_*\setminus p_* L$, index $*$ equals $\bullet$ or $\circ$.}

{\it Proof.} 
Let us denote the r.h.s. of (\ref{repL12}) as $L_{\bullet\circ}$ and show that it is an index $p_\circ$ sublatice of $L_\bullet$. 
Indeed, 
$$
L_{\bullet\circ}=Span_{\mathbb Z}\{p_\circ L_\bullet,~p_\bullet e_\circ\}
$$
and it remains to check that $p_\bullet e_\circ \not\in p_\circ L_\bullet$. If it were not the case, then $p_\bullet e_\circ\in$
$p_\bullet L\cap p_\circ L=$ $p_\bullet p_\circ L$ and $e_\circ\in p_\circ L$ contrary to our choice of $e_\circ$.
In the same fashion we check that $L_{\bullet\circ}$ is an index $p_\bullet$ sublattice of $L_\circ$.
We see that $L_{\bullet\circ}$ is a sublatice of the intersection $L_\bullet\cap L_\circ$.
Index of $L_\bullet\cap L_\circ$ in $ L$ is a multiple of both $p_\bullet$ and $p_\circ$, so it is at least $p_\bullet p_\circ$.
On the other hand $p_\bullet p_\circ=| L:L_{\bullet\circ}|=| L:L_\bullet\cap L_\circ||L_\bullet\cap L_\circ:L_{\bullet\circ}|$. Where from (\ref{repL12}) follows. ~~\bl

Combining Corollaries 1 and 2 we get the  following.
\begin{lemma}
Let $L_\bullet$ and $L_\circ$ be full rank sublattices of $L$ of prime indecies 
$p_\bullet$ and $p_\circ$ correspondingly and $L_{\bullet\circ}:=L_\bullet\cap L_\circ$. If
$L_\bullet\neq L_\circ$ then $|L_\bullet:L_{\bullet\circ}|=p_\circ$ and $|L_\circ:L_{\bullet\circ}|=p_\bullet$. 
Otherwise, if $L_\bullet= L_\circ$, then $|L_\bullet:L_{\bullet\circ}|= |L_\circ:L_{\bullet\circ}|=1$.
\label{ElementSquare}
\end{lemma}

Now we can list all prime index $p$ sublattices of $L$.
The factorset of any  sublattice (\ref{repL}) by its sublattice $pL$ consists of $p$ elements 
$\{je\}$, $j=0,\dots,p-1$, naturally included into the factorset $L/pL$ consisting of $p^2$ elements.
For different sublattices $L$, the factors $L/pL$ intersect only by the zero element  of $ L/pL$.
Therefore, there are exactly $(p^2-1)/(p-1)= p+1$ sublattices of prime index $p$
in $ L$. One can check that they are represented e.g. by the following transition matrices $Q$ for any fixed base in $ L$:
$$
\left(
\begin{array}{cc}
1&j\\0&p
\end{array}
\right),
\quad j=0,p-1,
\qquad
\left(
\begin{array}{cc}
p&0\\0&1
\end{array}
\right).
$$

\section{Composite index sublattices}
Let us fix an arbitrary lattice $ L$ and its full rank sublattices $ L_*$, $ L^*$.

For suitable bases in the lattices $ L$ and $ L_*$, the transition matrix $Q_*$ is diagonal (use Smith canonical form for integer matrix). Decomposing the elements of $Q_*$ into prime numbers we get a representation of the latter matrix as a product of integer matrices of prime determinants. Therefore we have the following chain of sublattices
$
 L:=L_0\supset L_1\supset L_2\dots\supset L_r=: L_*
$
of consecutive prime indecies $p_j:=|L_{j-1}:L_j|$. 
Same argument applied to the sublattice $ L^*$ gives us another filtration
$
 L:=L^0\supset L^1\supset L^2\dots\supset L^s=: L^*
$ with prime indecies $p^k:=|L^{k-1}:L^k|$.

We consider the sublattices $L_j^k:=L_j\cap L^k$ which naturally fill in the rectangular table

\be
\begin{array}{lllll}
 L^s\leftarrow & L_1^s\leftarrow & L_2^s\leftarrow & \dots \leftarrow &L_r^s= L_*\cap L^*:= L_*^*\\
\downarrow & \downarrow & \downarrow & \downarrow & \downarrow\\
\vdots&\vdots&\vdots& \vdots& \\
 L^2\leftarrow & L_1^2\leftarrow & \dots & \dots  \leftarrow &L_r^2\\
\downarrow & \downarrow & \downarrow & \downarrow & \downarrow\\
 L^1\leftarrow & L_1^1\leftarrow & L_2^1\leftarrow & \dots  \leftarrow & L_r^1\\
\downarrow & \downarrow & \downarrow & \downarrow & \downarrow\\
  L\leftarrow & L_1\leftarrow & L_2\leftarrow & \dots  \leftarrow &L_r\\
\end{array}
\label{Table}
\ee
where the arrows indicate the inclusions. Indeed, 
$$L_{j-1}^k\cap L_j^{k-1}:=(L_{j-1}\cap L^k)\cap(L_j\cap L^{k-1})=
(L_j\cap L_{j-1})\cap (L^k\cap L^{k-1})=L_j\cap L^k=:L_j^k.
$$

Applying lemma \ref{ElementSquare} consecutively to the elementary squares of the table (\ref{Table}) starting from the 
left-bottom one and moving  to the right along the lines of the table and upstairs along the columns we get the following

{\bf Corollory 3}  ~~~~~~~$|L_{j-1}^k:L_j^k|\in \{1,p_j\}$; ~~~~~~~~~~ $|L_j^{k-1}:L_j^k|\in \{1,p^k\}$.

\begin{theorem}
Any double decomposition (\ref{DD}) is the consequence of the relations of the same type 
with prime index sublattices $ L_\bullet$, $ L_\circ$.
\end{theorem}

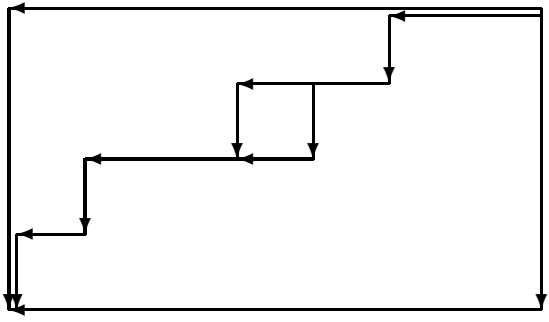
\begin{figure}
\begin{picture}(150, 40)
\put(-15,0){
\begin{picture}(120, 40) 
\thicklines
\put(120, 40){ \vector(-1,0){70} }
\put(50, 40){ \vector(0,-1){40} }
\put(120, 40){ \vector(0,-1){40} }
\put(120, 0){ \vector(-1,0){70} }

\put(120, 39){ \vector(-1,0){20} }
\put(100, 39){ \vector(0,-1){9} }
\put(100, 30){ \vector(-1,0){20} }
\put( 80, 30){ \vector(0,-1){10} }

\put( 90, 20){ \vector(-1,0){10} }
\put( 90, 30){ \vector(0,-1){10} }

\put( 80, 20){ \vector(-1,0){20} }
\put( 60, 20){ \vector(0,-1){10} }
\put( 60, 10){ \vector(-1,0){9} }
\put( 51, 10){ \vector(0,-1){10} }
\end{picture}}
\end{picture}
\caption{\small Deformation of paths on the table}
\label{Path}
\end{figure}

{\it Proof of Theorem 1.} Let us consider all possible paths coming from $ L_*^*$ to $ L$ 
along the arrows of the table (\ref{Table}). Each path corresponds to the filtration of 
the inital lattice $ L$ and therefore to the decomposition of the rational function
$R_{ L: L_*^*}(x)$ into prime compositional factors (including possibly identical elements).
The elementary change of the path caused by the alternative detour of the elementary square  in the table (see Fig. \ref{Path})
results in the change of two neighboring terms of the decomposition based on the double decomposition relation (\ref{DD})
$$
R_{L_j^k:L_{j+1}^k}\circ R_{L_{j+1}^k:L_{j+1}^{k+1}}=
R_{L_j^k:L_j^{k+1}}\circ R_{L_j^{k+1}:L_{j+1}^{k+1}}
$$
corresponding to prime index sublattices. The path coming along the top and left sides of the table 
may be converted to the path coming along the right and bottom sides by such elementary changes. ~~\bl

\parbox{9cm}
{\it
119991 Russia, Moscow GSP-1, ul. Gubkina 8,\\
Institute for Numerical Mathematics,\\
Russian Academy of Sciences\\[3mm]
{\tt gourmet@inm.ras.ru, ab.bogatyrev@gmail.com}}
\end{document}